\documentstyle[12pt]{article}

\newfont{\twelvemsb}{msbm10 scaled\magstep1}
\newfont{\eightmsb}{msbm8}
\newfont{\sixmsb}{msbm6}
\newfam\msbfam
\textfont\msbfam=\twelvemsb
\scriptfont\msbfam=\eightmsb
\scriptscriptfont\msbfam=\sixmsb
\catcode`\@=11
\def\Bbb{\ifmmode\let\next\Bbb@\else
  \def\next{\errmessage{Use \string\Bbb\space only in math mode}}\fi\next}
\def\Bbb@#1{{\Bbb@@{#1}}}
\def\Bbb@@#1{\fam\msbfam#1}
\newfont{\twelvegoth}{eufm10 scaled\magstep1}
\newfont{\tengoth}{eufm10}
\newfont{\eightgoth}{eufm8}
\newfam\gothfam
\textfont\gothfam=\twelvegoth
\scriptfont\gothfam=\eightgoth
\def\frak{\ifmmode\let\next\frak@\else
  \def\next{\errmessage{Use \string\frak\space only in math mode}}\fi\next}
\def\frak@#1{{\fam\gothfam{{#1}}}}
\catcode`@=12
\newcommand{\uqsmall}{U_q(\widehat{\mbox{\tengoth sl}}_2)}

\catcode`\@=11
\long\def\@makefntext#1{
\protect\noindent \hbox to 3.2pt {\hskip-.9pt
$^{{\ninerm\@thefnmark}}$\hfil}#1\hfill}                

\def\@makefnmark{\hbox to 0pt{$^{\@thefnmark}$\hss}}  

\def\ps@myheadings{\let\@mkboth\@gobbletwo
\def\@oddhead{\hbox{}
\rightmark\hfil\ninerm\thepage}
\def\@oddfoot{}\def\@evenhead{\ninerm\thepage\hfil
\leftmark\hbox{}}\def\@evenfoot{}
\def\sectionmark##1{}\def\subsectionmark##1{}}

\renewcommand{\thefootnote}{\fnsymbol{footnote}}

\newcounter{sectionc}\newcounter{subsectionc}\newcounter{subsubsectionc}
\renewcommand{\section}[1] {\vspace*{0.6cm}\addtocounter{sectionc}{1}
\setcounter{subsectionc}{0}\setcounter{subsubsectionc}{0}\noindent
        {\normalsize\bf\thesectionc. #1}\par\vspace*{0.4cm}}
\renewcommand{\subsection}[1] {\vspace*{0.6cm}\addtocounter{subsectionc}{1}
        \setcounter{subsubsectionc}{0}\noindent
        {\normalsize\it\thesectionc.\thesubsectionc. #1}\par\vspace*{0.4cm}}
\renewcommand{\subsubsection}[1]
{\vspace*{0.6cm}\addtocounter{subsubsectionc}{1}
     \noindent {\normalsize\rm\thesectionc.\thesubsectionc.\thesubsubsectionc.
        #1}\par\vspace*{0.4cm}}

\newcounter{appendixc}
\newcounter{subappendixc}[appendixc]
\newcounter{subsubappendixc}[subappendixc]

\renewcommand{\appendix}[1] {\vspace*{0.6cm}
        \refstepcounter{appendixc}
        \setcounter{figure}{0}
        \setcounter{table}{0}
        \setcounter{equation}{0}
        \renewcommand{\thefigure}{\Alph{appendixc}.\arabic{figure}}
        \renewcommand{\thetable}{\Alph{appendixc}.\arabic{table}}
        \renewcommand{\theappendixc}{\Alph{appendixc}}
        \renewcommand{\theequation}{\Alph{appendixc}.\arabic{equation}}
        \noindent{\bf Appendix \theappendixc #1}\par\vspace*{0.4cm}}

\def\abstracts#1{{
\centering{\begin{minipage}{12.2truecm}\footnotesize\baselineskip=12pt\noindent
        \centerline{\footnotesize ABSTRACT}\vspace*{0.3cm}
        \parindent=0pt #1
        \end{minipage}}\par}}

\renewenvironment{thebibliography}[1]
        {\begin{list}{\arabic{enumi}.}
        {\usecounter{enumi}\setlength{\parsep}{0pt}
\setlength{\leftmargin 1.25cm}{\rightmargin 0pt}
         \setlength{\itemsep}{0pt} \settowidth
        {\labelwidth}{#1.}}}{\end{list}}

\newcounter{itemlistc}
\newcounter{romanlistc}
\newcounter{alphlistc}
\newcounter{arabiclistc}

\newcommand{\fcaption}[1]{
        \refstepcounter{figure}
        \setbox\@tempboxa = \hbox{\footnotesize Fig.~\thefigure. #1}
        \ifdim \wd\@tempboxa > 6in
           {\begin{center}
        \parbox{6in}{\footnotesize\baselineskip=12pt Fig.~\thefigure. #1}
            \end{center}}
        \else
             {\begin{center}
             {\footnotesize Fig.~\thefigure. #1}
              \end{center}}
        \fi}

\newcommand{\tcaption}[1]{
        \refstepcounter{table}
        \setbox\@tempboxa = \hbox{\footnotesize Table~\thetable. #1}
        \ifdim \wd\@tempboxa > 6in
           {\begin{center}
        \parbox{6in}{\footnotesize\baselineskip=12pt Table~\thetable. #1}
            \end{center}}
        \else
             {\begin{center}
             {\footnotesize Table~\thetable. #1}
              \end{center}}
        \fi}

\def\@citex[#1]#2{\if@filesw\immediate\write\@auxout
        {\string\citation{#2}}\fi
\def\@citea{}\@cite{\@for\@citeb:=#2\do
        {\@citea\def\@citea{,}\@ifundefined
        {b@\@citeb}{{\bf ?}\@warning
        {Citation `\@citeb' on page \thepage \space undefined}}
        {\csname b@\@citeb\endcsname}}}{#1}}

\newif\if@cghi
\def\cite{\@cghitrue\@ifnextchar [{\@tempswatrue
        \@citex}{\@tempswafalse\@citex[]}}
\def\citelow{\@cghifalse\@ifnextchar [{\@tempswatrue
        \@citex}{\@tempswafalse\@citex[]}}
\def\@cite#1#2{{$\null^{#1}$\if@tempswa\typeout
        {IJCGA warning: optional citation argument
        ignored: `#2'} \fi}}

\catcode`@=12

 1
 1
 1

\font\ninerm=cmr9

\renewcommand{\theequation}{\thesectionc.\arabic{equation}}
\def\be{\begin{equation}}
\def\en{\end{equation}}
\def\bea{\begin{eqnarray}}
\def\ena{\end{eqnarray}}
\def\bean{\begin{eqnarray*}}
\def\enan{\end{eqnarray*}}
\def\lb#1{\label{eqn:#1}}
\def\rf#1{(\ref{eqn:#1})}
\def\F{{\cal F}}
\def\mod{{\rm mod}~}
\newcommand{\gsl}{{\frak {sl}}}
\newcommand{\slt}{\gsl_2}
\newcommand{\slh}{\widehat{\gsl}}
\def\slth{\slh_2}
\def\uq{U_q(\slth)}
\def\End{{\rm End}}
\def\ad{{\rm ad}\,}

\def\tPhi{\tilde\Phi}

\def\Z{{\Bbb Z}}
\def\C{{\Bbb C}}

\def\e{\varepsilon}
\def\vep{\varepsilon}
\def\z{\zeta}

\begin{document}

\baselineskip=15pt

\begin{flushright}
  RIMS-1016\\
  q-alg/9505009\\
  May 1995\\
  revised 26 September 1995
\end{flushright}
\vspace*{0.2cm}

\centerline{\normalsize\bf New Level-$0$ Action of $\uq$ on Level-$1$
  Modules\footnote{contribution to the Proceedings of the meeting
    `Statistical Mechanics and Quantum Field Theory', at the
    University of Southern California, May 16-21, 1994}}

\vspace*{0.3cm}
\centerline{\footnotesize Michio Jimbo}
\baselineskip=13pt
\centerline{\footnotesize\it Department of Mathematics, Faculty of Science,}
\baselineskip=12pt
\centerline{\footnotesize\it Kyoto University, Kyoto 606, Japan}
\centerline{\footnotesize E-mail: jimbo@kusm.kyoto-u.ac.jp}
\vspace*{0.3cm}

\centerline{\footnotesize Rinat Kedem}
\baselineskip=13pt
\centerline{\footnotesize\it Research Institute for Mathematical Sciences,}
\baselineskip=12pt
\centerline{\footnotesize\it Kyoto University, Kyoto 606, Japan}
\centerline{\footnotesize E-mail: rinat@kurims.kyoto-u.ac.jp}
\vspace*{0.3cm}

\centerline{\footnotesize Hitoshi Konno}
\baselineskip=13pt
\centerline{\footnotesize\it Yukawa Institute for Theoretical Physics,}
\baselineskip=12pt
\centerline{\footnotesize\it Kyoto University, Kyoto 606, Japan}
\centerline{\footnotesize E-mail: konno@yukawa.kyoto-u.ac.jp}
\vspace*{0.3cm}

\centerline{\footnotesize Tetsuji Miwa}
\baselineskip=13pt
\centerline{\footnotesize\it Research Institute for Mathematical Sciences,}
\baselineskip=12pt
\centerline{\footnotesize\it Kyoto University, Kyoto 606, Japan}
\centerline{\footnotesize E-mail: miwa@kurims.kyoto-u.ac.jp}

\vspace*{0.3cm}
\centerline{\footnotesize and}
\vspace*{0.3cm}

\centerline{\footnotesize Jens-Ulrik H. Petersen}
\baselineskip=13pt
\centerline{\footnotesize\it Research Institute for Mathematical Sciences,}
\baselineskip=12pt
\centerline{\footnotesize\it Kyoto University, Kyoto 606, Japan}
\centerline{\footnotesize E-mail: petersen@kurims.kyoto-u.ac.jp}

\vspace*{0.9cm}

\abstracts{We define a level-$0$ action of $\uqsmall$ on the sum of
  level-$1$ irreducible highest weight modules.  With the aid of the
  affine Hecke algebras, this action is realized on the basis created
  by the vertex operators. This is a $q$-analogue of the Yangian
  symmetry in conformal field theory.}

\normalsize\baselineskip=15pt
\setcounter{footnote}{0}
\renewcommand{\thefootnote}{\alph{footnote}}

\setcounter{equation}{0}
\section{Introduction}
At the USC meeting, a new symmetry structure in conformal field theory
(CFT), the Yangian symmetry of the spinon basis, was
introduced\cite{BLS94}. The aim of this paper is to investigate the
$q$-deformation of this mathematical structure. Before stating the
main achievement in this paper, let us briefly review the developments
in this subject.

Results suggesting a possible fermionic description of the CFT
spectrum were first obtained by Ref.\cite{SUNY} for certain critical
lattice models.  Closely related results appeared earlier in the
mathematical literature (see e.g., Ref.\cite{LP}).  Analysis of the
spectrum of these models from the Bethe ansatz yielded fermionic
expressions for the conformal partition functions (see
Refs.\cite{KKMM1,BONN,KNS,Mel,FeON,Bkv,BeMc,War,FoQu} for further
developments). This suggested a possible alternative description of
the states in CFT, the spinon basis of Refs.\cite{BPS,BLSa}.

The underlying algebraic structure, the Yangian symmetry, has been
suggested from quite a different direction.  The relevant subject was
the diagonalization of the Calogero--Sutherland model with spin degrees
of freedom\cite{Cal,Suth,BGHP} and the Haldane--Shastry
model\cite{Hal,Shas}.  Their spectrum is described by an asymptotic
Bethe ansatz\cite{Suth}.  However the multiplicity of states had never
been explained until the works\cite{HHTBP,BGHP}, which show that the
model has a symmetry generated by a Yangian\cite{Dri85}. The
multiplicities are identified with the dimensions of the highest
weight representations of the Yangian\cite{Dri88,ChPr90}.

\pagestyle{plain}

It turns out that the level-$1$ integrable highest weight module
(IHWM) of the affine Lie algebra $\slth$ has the same Yangian
structure as the space of states of these long range interaction
models at a special value of the coupling constant.  This fact was
established by constructing the new spinon basis in the level-$1$
IHWM\cite{BPS,BLSa} and considering the action of the Yangian on it.
The spinon basis is generated by the spin-1/2 primary field operator in the
level-$1$ $su(2)$ WZW model.  The character formula derived in the
spinon basis coincides with the fermionic representation. Recently,
the extension of spinon basis to the higher level case has also been
investigated\cite{BLSb}.

Now let us turn to our motivation in this work. In
Refs.\cite{FaTa81,FaTa84}, Faddeev and Takhtajan analysed the
structure of the space of the states in the XXX model in the infinite
volume limit. In Refs.\cite{DFJMN,JM}, the symmetry structure of the
off-critical XXZ model in the infinite volume limit was determined by
means of the representation theory of the quantum affine algebra
$\uq$.  The structure of the space of states was
understood as the level-$0$ representation on $\End_\C({\cal H})$
where ${\cal H}$ is the direct sum of the level-$1$ irreducible
highest weight $\uq$-modules.  The vacuum states are
identified with certain degree operators, and the excited states are
created by the type II vertex operator acting on the former.

This is similar to the Yangian structure we mentioned above, in the sense
that we have a level-$0$ action. However, the space
$\End_\C({\cal H})$ is much bigger than the space ${\cal H}$
itself, which goes to the IHWM in the conformal limit.
It is still unclear how the level-$0$ structure
on $\End_\C({\cal H})$ is related to the
Yangian structure on ${\cal H}$
in conformal field theory. As for this point, see the
recent paper by Nakayashiki and Yamada\cite{NaYa}.

The mathematical content of the Yangian structure was known since
Drinfeld and others\cite{Dri86b,Cher87}. Namely, there is a functorial
construction of Yangian representations from certain representations
of the degenerate affine Hecke algebra. In Ref.\cite{CP}, Chari and
Pressley gave a construction of level-$0$ $U_q(\slh_n)$-modules from
finite-dimensional representations of the affine Hecke algebra. In
this paper, we apply their construction to the level-$1$ $\uq$-module
${\cal H}$ and obtain a new level-$0$ action on it. For this purpose,
following Ref.\cite{BLSa}, we use the basis of ${\cal H}$ created by
the type I vertex operator. The construction of such a basis was
discussed in Ref.\cite{hwm} in a more general context corresponding to
the XYZ model.  We do not know if the `new' level-$0$ action of $\uq$ on
level-$1$ modules is related to some physical models, in a way the
Yangian structure is.

Before closing the introduction let us point out some technical details,
except for which our construction is a straightforward generalization of the
methods given in Refs.\cite{CP,BPS,BGHP,BLSa}.

(i) We are forced to consider infinite
sums of products of the components of vertex operators.
We will introduce suitable completion of the `$N$-spinon' space
to clarify the mathematical content.

(ii) Because of the fusion relation, which connects the `$N$-spinon'
sector to the `$(N-2)$-spinon' sector, we must check the compatibility
of the level-$0$ action in these different sectors. This is technically
non-trivial.

\setcounter{equation}{0}
\def\H{{\cal H}}
\def\V{{\cal V}}
\def\Ker{{\rm Ker}\,}
\def\ket#1{{|#1\rangle}}
\def\vaf{V_{\hbox{\scriptsize aff}}}
\def\op{\scriptstyle {\rm op}}
\def\uop{U^{\op}}
\def\ufin{U_q\left(\slt\right)\,}
\def\ufop{U_q^{\op}\left(\slt\right)\,}

\def\pf{\noindent{\it Proof.\quad}}
\def\qed{\hfil\fbox{}\medskip}
\newtheorem{prop}{Proposition}
\newtheorem{lem}[prop]{Lemma}

\section{Vertex operators}
The purpose of this section is to fix the notation concerning
level-$1$ modules and vertex operators.
Throughout this paper we fix a complex number $q$ such that $0<|q|<1$.

\subsection{The module $V_{\scriptstyle{\rm aff}}$}

Let  $U=\uq$  denote the quantum affine algebra with the standard generators
 $e_i,f_i,t_i=q^{h_i}$ ($i=0,1$) and $q^d$.
We shall identify $\ufin$ with the subalgebra of
$U$ generated by $e_1,f_1,t_1$.
We retain the coproduct $\Delta$ and the antipode $a$ in Ref.\cite{JM}, e.g.
\begin{eqnarray*}
&&\Delta(e_i)=e_i\otimes 1+t_i\otimes e_i,
\qquad
\Delta(f_i)=f_i\otimes t_i^{-1}+ 1\otimes f_i,
\\
&&\Delta(q^h)=q^h\otimes q^h
\qquad (h=h_i,d).
\end{eqnarray*}
In this paper we shall also use the opposite coproduct
 $\Delta^{\scriptstyle {\rm op}}=\sigma\circ\Delta$
($\sigma x\otimes y=y\otimes x$).
We let $\uop$ denote the algbera $U$ equipped with
the opposite coalgebra structure determined by $\Delta^{\op}$.

Let $V=\C v_+\oplus \C v_-$ be the irreducible two-dimensional module
over $\ufin$.
The affinization of $V$ is the $U$-module
\[
\vaf={\rm span}_\C \{ v_{\vep,n}\mid \vep=\pm,n\in \Z\,\}
\]
given as follows.
\begin{eqnarray*}
&&e_0 v_{+,n}=v_{-,n+1}, \qquad
f_0 v_{-,n}=v_{+,n-1},
\\
&&e_1 v_{-,n}=v_{+,n}, \qquad
f_1 v_{+,n}=v_{-,n},
\\
&&t_0^{-1} v_{\pm,n}=t_1 v_{\pm,n}=q^{\pm 1}v_{\pm,n},
\qquad
q^d v_{\pm,n}=q^n v_{\pm,n}.
\end{eqnarray*}
In terms of the generating series
\[
v_{\vep}(z)=\sum_{n}v_{\vep,n}z^{-n},
\]
the action of $U$ reads $e_0v_+(z)=zv_-(z)$, $e_1v_-(z)=v_+(z)$, and so on.

\subsection{Vertex operators}
Let  $V(\Lambda_i)=U\ket{i}$
be the integrable level-$1$ module of $U$
with highest weight vector $\ket{i}$ ($i=0,1$).
Set $\H=V(\Lambda_0)\oplus V(\Lambda_1)$.
We have the homogeneous gradation $\H=\bigoplus_{r\ge 0}\H_{-r}$
such that $\H_0=\C\ket{0}\oplus\C\ket{1}$.

For an element $f\in \End_\C\H$, the
(opposite) adjoint action of $x\in U$ is given by
\[
\ad\!^{\op} x.f=\sum x_{(2)}~f~ a^{-1}\bigl(x_{(1)}\bigr)
\]
where $\Delta(x)=\sum x_{(1)}\otimes x_{(2)}$.
For $f,g\in\End_\C\H$ and $u\in\H$ we have
\begin{eqnarray}
\ad\!^{\op} x.  f\circ g& =&\sum \ad\!^{\op} x_{(2)}.
f \circ  \ad\!^{\op} x_{(1)}.  g,
\label{eqn:adj1}\\
x. f(u) &=&\sum \ad\!^{\op} x_{(2)}.  f \left( x_{(1)}u\right).
\label{eqn:adj2}
\end{eqnarray}
With this action we regard $\End_\C \H$ as an $\uop$-module.

\begin{prop}
There exists a unique embedding of the $\uop$-module
\begin{equation}\label{eqn:embed}
\vaf\longrightarrow \End_\C\H,
\qquad
v_{\vep,n}\mapsto \tilde{\Phi}^*_{\vep,n}
\end{equation}
such that
\[
\tilde{\Phi}^*_{+,0}\ket{0}=\ket{1},
\quad
\tilde{\Phi}^*_{-,0}\ket{1}=\ket{0}.
\]
\end{prop}

The elements $\tilde{\Phi}^*_{\vep,n}\in\End_\C\H$ will be referred to
as components of the vertex operator of type I.\ To make contact with
the usual formulation, introduce the generating series corresponding
to $v_\vep(z)$
\[
\tilde{\Phi}^{*}_{\vep}(z)=\sum_n
\tilde{\Phi}^{*}_{\vep,n} z^{-n}.
\]
Then (\ref{eqn:embed}) is $\uop$-linear if and only if
(i) $\tilde{\Phi}^*_{\vep,n}\H_{r}\subset \H_{r+n}$ for all $n,r$,
and (ii) the following map is $U'$-linear:
\[
u\otimes v_\vep(z) \mapsto \tilde{\Phi}^*_{\vep}(z)u.
\qquad (u\in \H)
\]
Here $U'$ signifies the subalgebra generated by
 $e_i,f_i,t_i$ ($i=0,1$) with the original coproduct $\Delta$.
Hence $\tilde{\Phi}^*_\vep(z)$ has the same meaning as the
one employed in Ref.\cite{JM}.

For reasons of weights, $\tilde{\Phi}^*_{\vep,n}$ sends
$V(\Lambda_i)$ to $V(\Lambda_{1-i})$.
We will sometimes write the restriction to $V(\Lambda_i)$
as
\[
\tilde{\Phi}^{*(1-i,i)}_{\vep,n}~:~
V(\Lambda_i)\longrightarrow V(\Lambda_{1-i}).
\]

\subsection{Spanning vectors}
Set $\V_N=\vaf^{\otimes N}$, $\V=\bigoplus_{N\ge 0}\V_N$.
The embedding (\ref{eqn:embed}) gives rise to an $\uop$-linear map
\begin{eqnarray}
&& \rho~:~\V~~ \longrightarrow ~~\End_\C\H,
\label{eqn:rho}\\
&&v_{\vep_1,n_1}\otimes \cdots \otimes v_{\vep_N,n_N}
{}~~\mapsto~~
\tilde{\Phi}^*_{\vep_1,n_1}\cdots \tilde{\Phi}^*_{\vep_N,n_N}.
\nonumber~
\end{eqnarray}
Acting on the highest weight vector $\ket{0}$ we obtain
\begin{eqnarray}
&&\rho_0~:~\V~~ \longrightarrow ~~\H,
\label{eqn:rhoz}\\
&&v_{\vep_1,n_1}\otimes \cdots \otimes v_{\vep_N,n_N}
{}~~\mapsto~~
\tilde{\Phi}^*_{\vep_1,n_1}\cdots \tilde{\Phi}^*_{\vep_N,n_N}\ket{0}.
\label{eqn:spinon}
\end{eqnarray}

\begin{prop}
The map $\rho_0$ is $\ufop$-linear and surjective.
\end{prop}
\pf
The $\ufop$-linearity follows from (\ref{eqn:adj2}) and the fact that
$\ket{0}$ belongs to the trivial representation of $\ufin$.
Note that
the image of $\rho_0$ contains $\ket{0}, \ket{1}=\tilde{\Phi}^*_{+,0}\ket{0}$,
and hence $f_0\ket{0}=\tilde{\Phi}^*_{+,-1}\ket{1}$.
Using this and (\ref{eqn:adj2}) one checks readily that the image
is also invariant under the action of $U$.
The surjectivity follows from the irreducibility of $V(\Lambda_i)$.
\qed

The vectors (\ref{eqn:spinon}) thus constitute a spanning set of $\H$,
and we have an isomorphism of $\ufop$-modules
$\V/\Ker\rho_0 {\buildrel \sim \over \longrightarrow} \H$.
In section 3 we will determine $\Ker\rho_0$.

\subsection{Relations among vertex operators}
Let us summarize here the properties of vertex operators
which will be used later.
To state them we introduce the generating series
\begin{equation}
\varphi_{\vep_1,\cdots,\vep_N}(z_1,\cdots,z_N)
={\prod_{j=1}^N z_j^{(N-j-p_j+p_N)/2}\over\prod_{j<k}\eta(z_k/z_j)}
\tilde{\Phi}^{*(p_0,p_1)}_{\vep_1}(z_1)\cdots
\tilde{\Phi}^{*(p_{N-1},p_N)}_{\vep_N}(z_N).
\label{eqn:op}
\end{equation}
Here $p_j=0$ ($j\equiv N \mod 2$), $=1$  ($j\not \equiv N \mod 2$),
and we have set
\[
\eta(z)={(q^6z;q^4)_\infty\over(q^4z;q^4)_\infty},
\qquad
(z;p)_\infty=\prod_{n=0}^\infty(1-p^nz).
\]
Note that (\ref{eqn:op}) comprises only integral powers of $z_j$:
\[
\varphi_{\vep_1,\cdots,\vep_N}(z_1,\cdots,z_N)
=
\sum_{m_1,\ldots,m_N\in \Z}
\varphi_{\vep_1,\cdots,\vep_N;m_1,\cdots,m_N}z_1^{-m_1}\cdots z_N^{-m_N}.
\]
Each coefficient
$\varphi_{\vep_1,\cdots,\vep_N;m_1,\cdots,m_N}$
has a well-defined action on $\H$.

We use the following $R$-matrix $\tilde R(z)\in\End_\C V\otimes V$
regarded as a power series in $z$.
\begin{eqnarray*}
&&\tilde R(z)v_\vep\otimes v_\vep={z-q^2\over1-q^2z}v_\vep\otimes v_\vep,
\\
&&\tilde R(z)v_+\otimes v_-
={(1-q^2)z\over1-q^2z}v_+\otimes v_-+{q(z-1)\over1-q^2z}v_-\otimes v_+,
\\
&&\tilde R(z)v_-\otimes v_+
={q(z-1)\over1-q^2z}v_+\otimes v_-+{1-q^2\over1-q^2z}v_-\otimes v_+.
\end{eqnarray*}
In general,
for an element $A\in\End_\C V$ we set
\[
\pi_k(A)
\varphi_{\vep_1\cdots\vep_N,m_1\cdots m_N}
=\sum_{\vep'_k}
\varphi_{\vep_1\cdots\vep'_k\cdots\vep_N,m_1\cdots m_N}A_{\vep'_k\vep_k}
\]
where $Av_\vep=\sum_{\vep'}v_{\vep'}A_{\vep'\vep}$.
The Yang-Baxter equation for
$\tilde R_{j,k}(z)=(\pi_j\otimes\pi_k)\Bigl(\tilde R(z)\Bigr)$
then reads as follows.
\bean
&&\tilde R_{k,k+1}(z_2/z_1)
\tilde R_{k+1,k+2}(z_3/z_1)
\tilde R_{k,k+1}(z_3/z_2)\\
&&=\tilde R_{k+1,k+2}(z_3/z_2)
\tilde R_{k,k+1}(z_3/z_1)
\tilde R_{k+1,k+2}(z_2/z_1).
\enan

\begin{prop}\label{prop:op} The following relations hold:
\begin{description}
\item[Commutation Relation]
\begin{eqnarray}
&&\varphi_{\vep_1,\cdots,\vep_j,\vep_{j+1},\cdots,\vep_N}
(z_1,\cdots,z_{j+1},z_j,\cdots,z_N)
\label{eqn:phicom}\\
&&=\tilde{R}_{j,j+1}(z_{j+1}/z_j)
\varphi_{\vep_1,\cdots,\vep_j,\vep_{j+1},\cdots,\vep_N}
(z_1,\cdots,z_j,z_{j+1},\cdots,z_N).
\nonumber
\end{eqnarray}
\item[Fusion Relation]
\begin{eqnarray}
&&\varphi_{\vep_1,\cdots,\vep_N}(z_1,\cdots,z_N)\Bigr|_{z_{j+1}=q^{-2}z_j}
\nonumber\\
&&=(-q)^{N-j+(\vep_j-1)/2}\delta_{\vep_j+\vep_{j+1},0}
\prod_{i=1}^{j-1}(z_i-q^2z_j)\prod_{i=j+2}^{N}(q^{-2}z_j-q^2z_i)
\nonumber\\
&&\quad \times
\varphi_{\vep_1,\cdots,\vep_{j-1},\vep_{j+2},\cdots,\vep_N}
(z_1,\cdots,z_{j-1},z_{j+2},\cdots,z_N).\lb{FUS}
\end{eqnarray}
\item[Highest Weight Condition]
$\varphi_{\vep_1,\cdots,\vep_N}(z_1,\cdots,z_N)\ket{0}$ is a power series in
$z_N$.
\end{description}
\end{prop}

The first two are direct consequences of the corresponding
properties\cite{FR,DFJMN,JM} of $\tPhi^*_\vep(z)$, and the last one is
obvious.  We remark that, in the literature, the commutation relations
are written in the sense of (analytic continuation of) matrix
elements.  Thanks to the analyticity properties of type I vertex
operators, they are valid also as formal series (see the remark at the
end of A.4 in Ref.\cite{hwm}).  It is also possible to verify these
relations directly using bosonization (see Ref.\cite{JM}).

We end this section with a lemma, to be used later.
\begin{lem}\label{lem:max}
If $\max(m_1,\cdots,m_N)+r>0$, then
\begin{equation}
\varphi_{\vep_1\cdots\vep_N,m_1\cdots m_N}\H_r=0.
\label{eqn:ANN}
\end{equation}
\end{lem}
\pf
Suppose that $m_j=\max(m_1,\cdots,m_N)$.
If $j=N$ then $m_N+r>0$ and (\ref{eqn:ANN}) follows.
If $j<N$, then by applying (\ref{eqn:phicom}) for $k=j$ we can rewrite
the left hand side in the form
\begin{eqnarray*}
&&
\varphi_
{\vep_1\cdots\vep_j\vep_{j+1}\cdots\vep_N,m_1\cdots m_jm_{j+1}\cdots m_N}
\\
&&=\sum_{\vep'_j,\vep'_{j+1}\atop k\ge0}
\varphi_{\vep_1\cdots\vep'_j\vep'_{j+1}\cdots\vep_N,m_1\cdots,m_{j+1}-k,m_j+k,
\cdots m_N}
c^k_{\vep'_j\vep'_{j+1},\vep_j\vep_{j+1}}.
\end{eqnarray*}
Since
$\max(m_1,\cdots,m_{j+1}-k,m_j+k,\cdots,m_N)=m_j+k$,
the statement follows by induction.
\hfill \qed

\noindent{\sl Remark.}\quad
Taking $r=0$ we find that
$\varphi_{\vep_1,\cdots,\vep_N}(z_1,\cdots,z_N)\ket{0}$ is actually
a power series in $z_1,\cdots,z_N$.

\setcounter{equation}{0}
\def\tPhi{\tilde{\Phi}}
\def\Vh{\widehat{\V}}
\def\Fb{\overline{F}}
\def\N{{\cal N}}

\section{Basis of level-$1$ modules}

\subsection{Completion of $\V_N$}
In order to study the kernel of $\rho_0$, we deal with
a generating series corresponding to (\ref{eqn:op}),
whose coefficients are certain
infinite sums of $v_{\vep_1,n_1}\otimes \cdots \otimes v_{\vep_N,n_N}$.
For the precise formulation we introduce a completion of
$\V_N=\vaf^{\otimes N}$.

For each $r\in \Z$, let
\[
\V^{(r)}_N={\rm span}_\C
\{ v_{\vep_1,m_1}\otimes\cdots\otimes v_{\vep_N,m_N}
\mid m_1+\cdots+m_N=r\,\}.
\]
Setting
\[
D(m_1,\cdots,m_N)=\max(m_1+m_2+\cdots+m_N,
\cdots,m_{N-1}+m_N, m_N)
\]
we define a decreasing filtration
\begin{eqnarray*}
&&\V^{(r)}_N\supset \cdots \supset \V^{(r)}_N[l]
\supset \V^{(r)}_N[l+1] \supset \cdots,
\\
&&
\V^{(r)}_N[l]={\rm span}_\C
\{  v_{\vep_1,m_1}\otimes\cdots\otimes v_{\vep_N,m_N} \in \V^{(r)}_N
\mid D(m_1,\cdots,m_N)\ge l\,\}.
\end{eqnarray*}
Denote by $\Vh^{'(r)}_N$ the completion of $\V^{(r)}_N$
with respect to this filtration.
\[
\Vh^{'(r)}_N=\lim_{\longleftarrow\atop l}
\V^{(r)}_N/\V^{(r)}_N[l].
\]
We set $\Vh'_N=\oplus_{r\in\Z}\Vh^{'(r)}_N$, $\Vh'=\oplus_{N\geq 0}\Vh'_N$.
Since $\tPhi^*_{\vep,n}\H_s\subset \H_{s+n}$
and $\H_s=0$ for $s>0$, we have
$\rho\left(\V^{(r)}_N[l]\right)\H_s=0$ if $l+s>0$.
It follows that the map $\rho$, and hence $\rho_0$,
 extends to the completion:
\[
\widehat\rho'~:~\Vh'\longrightarrow\End_\C\H,
\qquad
\widehat\rho'_0~:~\Vh'\longrightarrow \H.
\]

\subsection{Generating series $F_{\vep_1,\cdots,\vep_N}(z_1,\cdots,z_N)$}

Let $p_j$ and $\eta(z)$ be as in (\ref{eqn:op}).
We define an analogous generating series whose coefficients are
in $\Vh'_N$:
\begin{eqnarray}
F_{\vep_1,\cdots,\vep_N}(z_1,\cdots,z_N)
&=&\frac{\prod_{j=1}^N z_j^{(N-j-p_j+p_N)/2}}{\prod_{j<k}\eta(z_k/z_j)}
v_{\vep_1}(z_1)\otimes\cdots\otimes v_{\vep_N}(z_N)
\nonumber\\
&=&
\sum_{m_1,\cdots,m_N\in\Z}
F_{\vep_1\cdots\vep_N,m_1\cdots m_N}
z_1^{-m_1}\cdots z_N^{-m_N}.
\label{eqn:gen}
\end{eqnarray}
If we set
\begin{eqnarray*}
\Z_{+,N}&=&\{\sum_{k=1}^{N-1}n_k(0,\cdots,
{\underbrace{-1,1}_{k,k+1}},\cdots,0)\mid
n_1,\cdots,n_{N-1}\in\Z_{\ge0}\},\\
\kappa&=&(\kappa_1,\cdots,\kappa_N), \qquad
\kappa_j=(N-j-p_j+p_N)/2,
\end{eqnarray*}
then the Laurent coefficients have the form
\[
F_{\vep_1\cdots\vep_N,m_1\cdots m_N}=
\sum_{(n_1,\cdots,n_N)\in\atop
(m_1,\cdots,m_N)+\kappa+\Z_{+,N}}
c_{n_1,\cdots,n_N}
v_{\vep_1,n_1}\otimes\cdots\otimes v_{\vep_N,n_N}
\]
for some $c_{n_1,\cdots,n_N}\in\C$.
Since
\[
D(m_1,\cdots,m_j-1,m_{j+1}+1,\cdots,m_N)
\ge D(m_1,\cdots,m_j,m_{j+1},\cdots,m_N),
\]
the sum $F_{\vep_1\cdots\vep_N,m_1\cdots m_N}$ converges to an element
of $\Vh'_N[l]$ with $l=D(m_1,\cdots,m_N)$.
Conversely
\[
v_{\vep_1,n_1}\otimes\cdots\otimes v_{\vep_N,n_N}=
\sum_{(m_1,\cdots,m_N)\atop\in
(n_1,\cdots,n_N)-\kappa+\Z_{+,N}}
\tilde{c}_{m_1,\cdots,m_N}
F_{\vep_1\cdots\vep_N,m_1\cdots m_N}
\]
for some $\tilde{c}_{m_1,\cdots,m_N}\in\C$.

\subsection{Second completion}

By the definition we have
\[
\widehat\rho'\left(F_{\vep_1\cdots\vep_N,m_1\cdots m_N}\right)
=
\varphi_{\vep_1\cdots\vep_N,m_1\cdots m_N}
\quad \in \End_\C\H.
\]
Therefore the commutation relation \rf{phicom} entails

\begin{prop}\label{prop:Frel}
The Laurent coefficients of the following belong to $\Ker \widehat\rho'$.
\begin{eqnarray}
&&F_{\vep_1,\cdots,\vep_j,\vep_{j+1},\cdots,\vep_N}
(z_1,\cdots,z_{j+1},z_j,\cdots,z_N)\nonumber\\
&&\quad -\tilde{R}_{j, j+1}(z_{j+1}/z_j)
F_{\vep_1,\cdots,\vep_j,\vep_{j+1},\cdots,\vep_N}
(z_1,\cdots,z_j,z_{j+1},\cdots,z_N).
\label{eqn:Fcom}
\end{eqnarray}
\end{prop}

Let us introduce the second filtration in $\Vh_N^{'(r)}$.
\begin{eqnarray*}
&&\Vh^{'(r)}_N\supset \cdots \supset \Vh^{'(r)}_N[[m]]
\supset \Vh^{'(r)}_N[[m+1]] \supset \cdots,
\\
&&\Vh^{'(r)}_N[[m]]={\rm cl.}
{\rm span}_\C \{ F_{\vep_1\cdots\vep_N,m_1\cdots m_N}\in \Vh^{'(r)}_N
\mid \max(m_1,\cdots,m_N)\ge m\,\}
\end{eqnarray*}
where cl stands for the closure in $\Vh'$.
We denote by $\Vh^{(r)}_N$ the completion of
$\Vh^{'(r)}_N$ with respect to this filtration
\[
\Vh_N^{(r)}=
\lim_{\longleftarrow\atop m}\Vh^{'(r)}_N/\Vh^{'(r)}_N[[m]],
\]
and set $\Vh_N=\oplus_{r\in\Z}\Vh^{(r)}_N$, $\Vh=\oplus_{N\geq 0}\Vh_N$.

In view of Lemma \ref{lem:max} and Proposition \ref{prop:Frel}, we see
that the maps $\widehat\rho'$ and $\widehat\rho'_0$ extend further to $\Vh$:
\[
\widehat\rho~:~\Vh\longrightarrow\End_\C\H,
\qquad
\widehat\rho_0~:~\Vh\longrightarrow \H.
\]

After preparing the second completion we can state (see \rf{FUS})
\begin{prop}
The Laurent coefficients of the following belong to $\Ker \widehat\rho$.
\begin{eqnarray}
&& F_{\vep_1,\cdots,\vep_N}(z_1,\cdots,z_N)\Bigr|_{z_{j+1}=q^{-2}z_j}
\nonumber\\
&&\quad -(-q)^{N-j+(\vep_j-1)/2}\delta_{\vep_j+\vep_{j+1},0}
\prod_{i=1}^{j-1}(z_i-q^2z_j)\prod_{i=j+2}^{N}(q^{-2}z_j-q^2z_i)
\nonumber\\
&&\qquad \times
F_{\vep_1,\cdots,\vep_{j-1},\vep_{j+2},\cdots,\vep_N}
(z_1,\cdots,z_{j-1},z_{j+2},\cdots,z_N).
\label{eqn:Ffus}
\end{eqnarray}
\end{prop}

\subsection{Kernel of $\widehat\rho_0$}

{}From the remark at the end of 2.4 we see that
$\Ker\widehat\rho_0$ contains
\begin{equation}
F_{\vep_1\cdots\vep_N,m_1\cdots m_N}\quad \hbox{with $m_j>0$ for some $j$.}
\label{eqn:Fhwt}
\end{equation}
The Laurent coefficients of (\ref{eqn:Fcom}) and (\ref{eqn:Ffus})
reduce to finite sums modulo
$\Ker\widehat\rho_0$,
and therefore
$\Ker\rho_0$ contains (\ref{eqn:Fcom}) and (\ref{eqn:Ffus})
in reduced form. Thus, it follows that
\[
\widehat\V/\Ker\widehat\rho_0\simeq
\widehat\V'/\Ker\widehat\rho'_0\simeq
\V/\Ker\rho_0\simeq\H.
\]

Let $\N\subset \Vh$ be (the closure of) the
span of elements (\ref{eqn:Fcom}), (\ref{eqn:Ffus}) and (\ref{eqn:Fhwt}).

\begin{prop}
\[
\widehat\V/\N\simeq\H.
\]
\end{prop}

This can be verified as follows.  As we will show in Proposition 8, the general
expression $v_{\vep_1,n_1}\otimes\cdots\otimes v_{\vep_N, n_N}$ can be
reduced modulo $\N$ to a certain normal ordered form.  We then count
the number of such normal ordered expressions to find that the
character of $\Vh/\N$ is dominated by that of $\H$.  Since
$\widehat\rho_0:\Vh\rightarrow \H$ is surjective, this proves the
proposition.  In fact this working has already been discussed in the
paper\cite{hwm}.  (Note that the properties of the elliptic algebra
assumed in Ref.\cite{hwm} are all valid for quantum affine algebras.)
Hence it suffices to show that the `normal ordering rules' of
Ref.\cite{hwm} are consequences of (\ref{eqn:Fcom}) and
(\ref{eqn:Ffus}).

Define
\begin{equation}
\Fb_{\vep_1,\cdots,\vep_N}(\z_1,\cdots,\z_N)
={\prod_{j=1}^N\z_j^{(1+\vep_j)/2}\over
\prod_{j=1}^Nz_j^{N-j}\prod_{j<k}(1-q^2z_k/z_j)}
F_{-\vep_1\cdots-\vep_N}(z_1,\cdots,z_N)
\label{eqn:g}
\end{equation}
where $z_j=\z_j^2$.
Note that the coefficients of (\ref{eqn:g}) are well--defined in $\Vh$.

\begin{prop}
The following belong to $\Ker\widehat\rho$.
\bea
&&\Fb_{\vep_1\cdots{\underbrace{\scriptstyle\vep\vep}_{k,k+1}}\cdots\vep_N}
(\z_1,\cdots,\z_{k+1},\z_k,\cdots,\z_N)
\nonumber\\
&&\qquad -
\Fb_{\vep_1\cdots{\underbrace{\scriptstyle\vep\vep}_{k,k+1}}\cdots\vep_N}
(\z_1,\cdots,\z_k,\z_{k+1},\cdots,\z_N),
\lb{A}\\
&&(\z_{k+1}+q\z_k)\Bigl(
\Fb_{\vep_1\cdots{\underbrace{\scriptstyle +-}_{k,k+1}}\cdots\vep_N}
(\z_1,\cdots,\z_{k+1},\z_k,\cdots,\z_N)
\nonumber\\
&&\qquad +\Fb_{\vep_1\cdots{\underbrace{\scriptstyle -+}_{k,k+1}}\cdots\vep_N}
(\z_1,\cdots,\z_{k+1},\z_k,\cdots,\z_N)\Bigr)
\nonumber\\
&&-
(\z_k+q\z_{k+1})\Bigl(
\Fb_{\vep_1\cdots{\underbrace{\scriptstyle +-}_{k,k+1}}\cdots\vep_N}
(\z_1,\cdots,\z_k,\z_{k+1},\cdots,\z_N)\nonumber\\
&&\qquad +\Fb_{\vep_1\cdots{\underbrace{\scriptstyle -+}_{k,k+1}}\cdots\vep_N}
(\z_1,\cdots,\z_k,\z_{k+1},\cdots,\z_N)\Bigr).
\lb{B}
\ena
\end{prop}
These relations are precisely the normal ordering rules in Ref.\cite{hwm}.
\smallskip

\pf
Let us write $A\sim B$ to mean that $A-B\in \Ker\widehat\rho$.
Then \rf{A} is shown as follows.
\bean
&&\Fb_{\vep_1\cdots{\underbrace{\scriptstyle\vep\vep}_{k,k+1}}\cdots\vep_N}
(\z_1,\cdots,\z_{k+1},\z_k,\cdots,\z_N)\\
&&
={\prod_{j=1}^N\z_j^{(1+\vep_j)/2}\over
\prod_{j=1}^Nz_j^{N-j}\prod_{j<k}(1-q^2z_k/z_j)}
{z_k(1-q^2z_{k+1}/z_k)\over z_{k+1}(1-q^2z_k/z_{k+1})}\\
&&\qquad \quad \times
F_{-\vep_1\cdots\underbrace{\scriptstyle -\vep,-\vep}_{k,k+1}\cdots-\vep_N}
(z_1,\cdots,z_{k+1},z_k,\cdots,z_N)\\
&&\sim
{\prod_{j=1}^N\z_j^{(1+\vep_j)/2}\over
\prod_{j=1}^Nz_j^{N-j}\prod_{j<k}(1-q^2z_k/z_j)}\\
&&\qquad \quad \times
F_{-\vep_1\cdots\underbrace{\scriptstyle -\vep,-\vep}_{k,k+1}\cdots-\vep_N}
(z_1,\cdots,z_k,z_{k+1},\cdots,z_N)\\
&&=\Fb_{\vep_1\cdots{\underbrace{\scriptstyle \vep\vep}_{k,k+1}}\cdots\vep_N}
(\z_1,\cdots,\z_k,\z_{k+1},\cdots,\z_N).
\enan
Similarly, we have
\bean
&&\Fb_{\vep_1\cdots{\underbrace{\scriptstyle \pm\mp}_{k,k+1}}\cdots\vep_N}
(\z_1,\cdots,\z_{k+1},\z_k,\cdots,\z_N)\\
&&\sim
{(1-q^2)\z_k/\z_{k+1}\over1-q^2z_k/z_{k+1}}
\Fb_{\vep_1\cdots{\underbrace{\scriptstyle \mp\pm}_{k,k+1}}\cdots\vep_N}
(\z_k,\cdots,\z_k,\z_{k+1},\cdots,\z_N)\\
&&+{q(1-z_k/z_{k+1})\over1-q^2z_k/z_{k+1}}
\Fb_{\vep_1\cdots{\underbrace{\scriptstyle \pm\mp}_{k,k+1}}\cdots\vep_N}
(\z_k,\cdots,\z_k,\z_{k+1},\cdots,\z_N).
\enan
The statement \rf{B} follows from this.
\qed

\subsection{Hecke algebra}

In section 4 we will define a `new' level-$0$ action of $U$ on the level-$1$
module $\H$.
For this purpose we need to rewrite the commutation relation (\ref{eqn:Fcom})
in the language of Hecke algebras.

Define $S\in\End_\C V\otimes V$ by
\begin{eqnarray*}
&&Sv_\vep\otimes v_\vep=-q^{-1}v_\vep\otimes v_\vep,\\
&&Sv_+\otimes v_-=(q-q^{-1})v_+\otimes v_--v_-\otimes v_+,\\
&&Sv_-\otimes v_+=-v_+\otimes v_-.
\end{eqnarray*}
The operators $S_{j,j+1}\in\End_\C V^{\otimes N}$ ($j=1,\cdots,N-1$)
satisfy the Hecke algebra relation
\begin{eqnarray}
&&S_{j,j+1}-S_{j,j+1}^{-1}=q-q^{-1},
\label{eqn:Hecke1}\\
&&S_{j,j+1}S_{k,k+1}=S_{k,k+1}S_{j,j+1}
\qquad (|j-k|>1),
\label{eqn:Hecke2}\\
&&S_{j,j+1}S_{j+1,j+2}S_{j,j+1}=S_{j+1,j+2}S_{j,j+1}S_{j+1,j+2}.
\label{eqn:Hecke3}
\end{eqnarray}
The $R$-matrix $\tilde R(z)$ can be written as
\begin{equation}
\tilde R(z)={Sz-S^{-1}\over qz-q^{-1}}.
\label{eqn:RS}
\end{equation}

Substituting (\ref{eqn:RS}) in (\ref{eqn:Fcom}) we obtain
formally the relations of the form
\begin{equation}
(S_{j,j+1}-G_{j,j+1})F_{\vep_1,\cdots,\vep_N}(z_1,\cdots,z_N)\sim0,
\label{eqn:SG}
\end{equation}
where
\begin{equation}
G_{j,k}^{\pm 1}={q^{-1}z_j-qz_k\over z_j-z_k}(K_{j,k}-1)+q^{\pm 1}
\label{eqn:G}
\end{equation}
and $K_{j,k}$ signifies the exchange of variables $z_j$ and $z_k$.
The $G_{j,j+1}$ ($j=1,\cdots,N-1$) satisfy
the same relations (\ref{eqn:Hecke1})--(\ref{eqn:Hecke3}) as do $S_{j,j+1}$.

The operators $G_{j,k}$ and $K_{j,k}$ are acting on formal Laurent series
$f(z_1,\cdots,z_N)$ in $z_1,\cdots,z_N$, and
preserve $\C[z_1,z_1^{-1},\cdots,z_N,z_N^{-1}]$.
However, even if the coefficients of $f(z_1,\cdots,z_N)$
all belong to $\Vh$,
those of $G_{j,k}f(z_1,\cdots,z_N)$ no longer do so in general.
In order to make sense of (\ref{eqn:SG}) we need
some consideration given below.

\subsection{Local series}

Consider a series in $(z_1,\cdots,z_N)$ with coefficients in $\Vh_N$,
\begin{equation}
f(z_1,\cdots,z_N)=
\sum_{m_1,\cdots,m_N\in\Z}f_{m_1,\cdots,m_N}z_1^{-m_1}\cdots z_N^{-m_N},
\label{eqn:local}
\end{equation}
where $f_{m_1,\cdots,m_N}\in\Vh_N$.
We call $f(z_1,\cdots,z_N)$ {\it local} and homogeneous
if there exists an $m\in \Z$ such that
\be
f_{m_1,\cdots,m_N}\in\Vh^{(m_1+\cdots m_N)}_N[[m+\max(m_1,\cdots,m_N)]],
\en
for all $m_1,\cdots,m_N$.
In general, we call $f(z_1,\cdots,z_N)$ local if it
can be written as a finite sum of local homogeneous series.
The generating series $F_{\vep_1,\cdots,\vep_N}(z_1,\cdots,z_N)$
(\ref{eqn:gen}) is local in this sense.

An important property of local series is the following:

\begin{lem}
Consider an arbitrary Laurent series homogeneous in $z_1,\cdots,z_N$:
\[
c(z_1,\cdots,z_N)=
\sum_{k_1,\cdots,k_N\in\Z\atop k_1+\cdots+k_N=l}
c_{k_1,\cdots,k_N}z_1^{k_1}\cdots z_N^{k_N}.
\]
If $f(z_1,\cdots,z_N)$ is local, then
the product $c(z_1,\cdots,z_N)f(z_1,\cdots,z_N)$ is well-defined
in $\Vh$.
\end{lem}

\pf
Expand the product
\begin{eqnarray*}
&&c(z_1,\cdots,z_N)f(z_1,\cdots,z_N)\\
&&=\sum_{m_1,\cdots,m_N}\left(
\sum_{k_1,\cdots,k_N\atop k_1+\cdots +k_N=l}
c_{k_1,\cdots,k_N}f_{m_1+k_1\cdots m_N+k_N}
\right)
z_1^{-m_1}\cdots z_N^{-m_N}.
\end{eqnarray*}
We must show that each coefficient
\[
f'_{m_1\cdots m_N}=
\sum_{k_1,\cdots,k_N\atop k_1+\cdots +k_N=l}
c_{k_1,\cdots,k_N}f_{m_1+k_1\cdots m_N+k_N}
\]
is convergent in $\Vh_N$.
This holds since there are only finitely many $k_j$'s
such that $\max(m_1+k_1,\cdots,m_N+k_N)$ is bounded from above and
$k_1+\cdots+k_N=l$.
\qed

In particular, we can multiply any homogeneous element in the ring of
Laurent series that are convergent on the unit circles
$|z_1|=\cdots=|z_N|=1$.
We note however that the product
$c(z_1,\cdots,z_N)f(z_1,\cdots,z_N)$ is not necessarily local.

\begin{lem}
If $f(z_1,\cdots,z_N)$ is local,
then $G_{j,k}f(z_1,\cdots,z_N)$ is well-defined and local.
\end{lem}
\pf
Applying $G_{j,k}$ to a monomial we obtain an expression
\[
G_{j,k}z_j^{-m_j}z_k^{-m_k} = \sum c^{m_jm_k}_{n_jn_k}z_j^{-n_j}z_k^{-n_k}
\]
where the sum is taken over $n_j,n_k$ such that
$n_j+n_k=m_j+m_k$ and $\max(n_j,n_k)\le \max(m_j,m_k)$.
Therefore if (\ref{eqn:local}) is  local and homogeneous, then
the coefficients of $G_{j,k}f(z_1,\cdots,z_N)$ have the form
\[
f'_{n_1\cdots n_N}= \sum_{m_j,m_k}
c^{m_jm_k}_{n_jn_k}f_{n_1,\cdots,m_j,\cdots,m_k,\cdots, n_N}
\]
where $m_j+m_k=n_j+n_k$ and $\max(m_j,m_k)\ge \max(n_j,n_k)$.
The lemma follows from this observation.
\qed

\subsection{The `$S=G$' relations}

Let us return to the relations (\ref{eqn:SG}).

\begin{prop}
The relations (\ref{eqn:Fcom}) are equivalent to
\begin{equation}
(S_{j,j+1}-G_{j,j+1})F_{\vep_1,\cdots,\vep_N}(z_1,\cdots,z_N)\sim0.
\lb{HEC}
\end{equation}
\end{prop}

\pf
Write (\ref{eqn:Fcom}) as
\begin{eqnarray}
&&F_{\vep_1,\cdots,\vep_j,\vep_{j+1},\cdots,\vep_N}
(z_1,\cdots,z_{j+1},z_j,\cdots,z_N)
\nonumber\\
&&\sim \tilde{R}_{j,j+1}(z_{j+1}/z_j)
F_{\vep_1,\cdots,\vep_j,\vep_{j+1},\cdots,\vep_N}
(z_1,\cdots,z_j,z_{j+1},\cdots,z_N).
\label{eqn:AAA}
\end{eqnarray}
Multiplying $qz_{j+1}-q^{-1}z_j$ we obtain
\begin{eqnarray}
&&(qz_{j+1}-q^{-1}z_j)F_{\vep_1,\cdots,\vep_j,\vep_{j+1},\cdots,\vep_N}
(z_1,\cdots,z_{j+1},z_j,\cdots,z_N)\nonumber\\
&&\sim(S_{j,j+1}z_{j+1}-S_{j,j+1}^{-1}z_j)
F_{\vep_1,\cdots,\vep_j,\vep_{j+1},\cdots,\vep_N}
(z_1,\cdots,z_j,z_{j+1},\cdots,z_N).
\nonumber\\
\lb{MID}
\end{eqnarray}
Conversely, by multiplying the Laurent series
$1/(qz_{j+1}-q^{-1}z_j)$, which is convergent for
$|z_j|=|z_{j+1}|=1$, we get (\ref{eqn:AAA}) from \rf{MID}.
Thus (\ref{eqn:AAA}) and \rf{MID} are equivalent.

The relation \rf{MID} is equivalent to
\begin{equation}
(z_{j+1}-z_j)\Bigl((S_{j,j+1}-G_{j,j+1})
F_{\vep_1,\cdots,\vep_N}(z_1,\cdots,z_N)\Bigr)\sim0.\lb{ZMUL}
\end{equation}
Set
\begin{equation}
(S_{j,j+1}-G_{j,j+1})F_{\vep_1,\cdots,\vep_N}(z_1,\cdots,z_N)
=\sum_{m_1,\cdots,m_N\in\Z}h_{m_1,\cdots,m_N}z_1^{-m_1}\cdots z_N^{-m_N}.
\lb{LOC}
\end{equation}
Then \rf{LOC} is local, and we have
\[
h_{m_1,\cdots,m_j,m_{j+1},\cdots,m_N}
\sim h_{m_1,\cdots,m_j-1,m_{j+1}+1,\cdots,m_N}.
\]
Using this repeatedly, we find $h_{m_1,\cdots,m_N}\sim0$.
\qed

\setcounter{equation}{0}
\def\P{{\cal P}}
\def\N{{\cal N}}
\def\End{\hbox{\rm End}}
\def\H{{\cal H}}
\def\F{\hat{\cal V}}
\def\V{{\cal V}}

\def\jj{\widehat{\scriptscriptstyle j,j+1}}
\def\res{\Big|_{z_N=q^{-2}z_{N-1}}}
\def\g#1#2{G^{-1}_{N-#1,N-#2}}
\def\es#1#2{S_{N-#1,N-#2}}
\def\oc{{\overline C}}

\newtheorem{thm}[prop]{Theorem}

\section{New action}

\subsection{Affine Hecke algebras}

The action of $\ufop$ on each $\V_N=\vaf^{\otimes N}$ extends
naturally to the completion $\hat\V$.  The subspace ${\cal
  N}\subset\hat\V$, given as the span of \rf{HEC}, (\ref{eqn:Ffus})
and (\ref{eqn:Fhwt}), is invariant with respect to this action.  This
can be seen explicitly by noting that $[S, \Delta^{\op}(x)]=0$
($\forall x\in\ufin$) and that the vector
\begin{equation}
v_+\otimes v_--q^{-1}v_-\otimes v_+
\label{eqn:sing}
\end{equation}
belongs to the trivial representation of $\ufin$.
Hence $\H\simeq \Vh/\N$ is an $\ufin$-module.
We will extend this to the level-$0$ action of the quantum
affine algebra $U$.

Recall that, given complex numbers $a_j\neq 0$ ($j=1,\cdots ,N$),
one can extend the $\ufin$-module $\V_N$ to
the evaluation module over $U$ (denoted by $\pi_{a_1,\cdots,a_N}$)
\begin{eqnarray}
&&\pi_{a_1,\cdots,a_N}(e_0)=\sum_{j=1}^Na_j\pi_j(f_1)
\pi_{j+1}(t_1^{-1})\cdots
\pi_N(t_1^{-1}),
\label{eqn:newe}\\
&&\pi_{a_1,\cdots,a_N}(f_0)=\sum_{j=1}^N a_j^{-1}\pi_1(t_1)
\cdots \pi_{j-1}(t_1)
\pi_j(e_1),
\label{eqn:newf}\\
&&\pi_{a_1,\cdots,a_N}(t_0)=
\pi_{1}(t_1^{-1})\cdots
\pi_N(t_1^{-1}).
\label{eqn:newt}
\end{eqnarray}
However this action does not descend to $\H$ because it violates the
condition $e_0\N\subset \N$ and $f_0\N\subset \N$.
The remedy is to replace the numbers $a_j$ by suitable
commuting operators, to be given below.

Following Ref.\cite{BGHP} let us introduce the operators
\begin{eqnarray}
Y_j&=&(G_{j,j+1}^{-1}K_{j,j+1})\cdots(G_{j,N}^{-1}K_{j,N})
p^{\vartheta_j}(K_{1,j}G_{1,j})\cdots(G_{j-1,j}K_{j-1,j})
\nonumber\\
&=&G_{j,j+1}^{-1}\cdots G_{N-1,N}^{-1}Z
G_{1,2}\cdots G_{j-1,j},
\lb{GENY}
\end{eqnarray}
where $Z=K_{1,2}K_{1,3}\cdots K_{1,N}p^{\vartheta_1}$
and $p^{\vartheta_j}$ denotes the scale operator
\[
p^{\vartheta_j}f(z_1,\cdots,z_N)=f(z_1,\cdots,pz_j,\cdots,z_N).
\]
At this stage the parameter $p$ is arbitrary, but
we will make a specific choice of it later on.

The operators $G_{j,j+1}$ ($j=1,\cdots, N-1$) and $Y_j$ ($j=1,\cdots,
N$) are known to satisfy the relations for the affine Hecke algebra
$\hat{H}_N$.  Namely we have, in addition to
(\ref{eqn:Hecke1})--(\ref{eqn:Hecke3}),
\begin{eqnarray*}
&&Y_jY_k=Y_kY_j, \\
&&G_{j,j+1} Y_j G_{j,j+1}= Y_{j+1},\\
&&[G_{j,j+1}, Y_k]=0, \qquad (j, j+1\not=k).
\end{eqnarray*}

In what follows,
for an operator $X\in\End_\C\C[z_1,z_1^{-1},\cdots,z_N,z_N^{-1}]$
we define $\hat X\in\End_\C\F_N$ by setting
\[
\sum \left(\hat X F_{\e_1\cdots\e_N;m_1,\cdots,m_N}\right)
z_1^{-m_1}\cdots z_N^{-m_N}
=XF_{\e_1\cdots\e_N}(z_1,\cdots,z_N).
\]
Notice that $\widehat{XY}=\hat{Y}\hat{X}$.

We now define $\pi^{(N)}(e_0)$, $\pi^{(N)}(f_0)$ by the right hand side
of (\ref{eqn:newe}) and (\ref{eqn:newf}) respectively, wherein we set
\begin{equation}
a_j=q^{N-1}\hat{Y}_j^{-1}.
\label{eqn:aj}
\end{equation}

The following says that they leave the subspace
$\N$ invariant.
\begin{prop}\label{prop:ess} Let $p=q^4$.
For $x=e_0$ or $f_0$, we have
\bea
&&\pi^{(N)}(x)\Bigl(S_{j,j+1}-G_{j,j+1}\Bigr)
F_{\e_1\cdots\e_N}(z_1,\cdots,z_N)\sim 0,
\lb{rhosg}\\
&&\left(\pi^{(N)}(x)
  F_{\e_1\cdots\e_N}(z_1,\cdots, z_j, z_{j+1},\cdots,z_N)\right)
\Bigl|_{z_{j+1}=q^{-2}z_j}
\nonumber\\
&&\quad \qquad \sim
(-q)^{N-j+(\vep_j-1)/2}\delta_{\vep_j+\vep_{j+1},0}
\prod_{i=1}^{j-1}(z_i-q^2z_j)\prod_{i=j+2}^{N}(q^{-2}z_j-q^2z_i)
\nonumber\\
&&\times \pi^{(N-2)}(x)
F_{\e_1\cdots\e_N}(z_1,\cdots,z_{j-1},z_{j+1},\cdots,z_N),
\lb{rhof}\\
&&\pi^{(N)}(x)F_{\e_1\cdots\e_N,m_1,\cdots,m_N}
\hbox{ \rm belongs to the kernel of $\widehat\rho_0$}\nonumber\\
&&\qquad\hbox{\rm if $m_j>0$ for some $j$}.\lb{HWCR}
\ena
\end{prop}

This proposition gives us

\begin{thm}
The formulas (\ref{eqn:newe})--(\ref{eqn:newt}) with the choice
(\ref{eqn:aj}), \rf{GENY} and $p=q^4$
defines a level-$0$ action of $U$ on $\H$.
\end{thm}

The rest of the text is devoted to the proof of Proposition \ref{prop:ess}.
The statement \rf{HWCR} is valid because the operator $Y_j$ keeps
the spce $\C[z_1^{-1},\cdots,z_N^{-1}]$ invariant and preserves the degree.

\subsection{Proof of \rf{rhosg}}

For this the parameter $p$ can be arbitrary.
The way of verification is the same as in Ref.\cite{CP}.
Consider the case $x=e_0$ and set
\begin{eqnarray*}
&&e^{(k)}_0=q^{N-1}\hat{Y}_k^{-1}f^{(k)}, \\
&&f^{(k)}=\pi_k(f_1)
\pi_{k+1}(t_1^{-1})\cdots
\pi_N(t_1^{-1}).
\end{eqnarray*}
Clearly $e_0^{(k)}$ commutes with $G_{j,j+1}$ and
$S_{j,j+1}$ if $k\neq j,j+1$.
Using
\be
S(f_1\otimes t_1^{-1})=(1\otimes f_1)S,
\lb{SFREL}
\en
and
$[\hat{Y}_k^{-1}, S_{j,j+1}]=0$, we have
\begin{eqnarray*}
e_0^{(j)}S_{j,j+1}&=&e_0^{(j)}(S_{j,j+1}^{-1}+q-q^{-1})\nonumber \\
&=&S_{j,j+1}^{-1}q^{N-1}\hat{Y}_j^{-1}f^{(j+1)}+(q-q^{-1})e_0^{(j)},
\\
e_0^{(j+1)}S_{j,j+1}&=&
S_{j,j+1}q^{N-1}\hat{Y}_{j+1}^{-1}f^{(j)}
\end{eqnarray*}
On the other hand, from
$\hat{G}_{j,j+1}\hat{Y}_{j}\hat{G}_{j,j+1}=\hat{Y}_{j+1}$ and
$[\hat{G}_{j,j+1},f^{(k)}]=0$, we find
\begin{eqnarray*}
e_0^{(j)}\hat{G}_{j,j+1}&=&e_0^{(j)}(\hat{G}_{j,j+1}^{-1}+q-q^{-1})\nonumber\\
&=&\hat{G}_{j,j+1}q^{N-1}\hat{Y}_{j+1}^{-1}f^{(j)}+(q-q^{-1})e_0^{(j)},
\\
e_0^{(j+1)}\hat{G}_{j,j+1}&=&
\hat{G}_{j,j+1}^{-1}q^{N-1}\hat{Y}_{j}^{-1}f^{(j+1)}.
\end{eqnarray*}
Therefore, we have
\begin{eqnarray*}
&&(e_0^{(j)}+e_0^{(j+1)})
\Bigl(S_{j,j+1}-\hat{G}_{j,j+1}\Bigr)
\\
&=&
\Bigl(S_{j,j+1}^{-1}-\hat{G}_{j,j+1}^{-1}\Bigr)
q^{N-1}\hat{Y}_{j}^{-1}f^{(j+1)}
\\
&&\qquad\qquad\qquad\qquad +\Bigl(S_{j,j+1}-\hat{G}_{j,j+1}\Bigr)
q^{N-1}\hat{Y}_{j+1}^{-1}f^{(j)}
\\
&=& (S_{j,j+1}-\hat{G}_{j,j+1}) q^{N-1}
  (\hat{Y}_{j}^{-1}f^{(j+1)}+ \hat{Y}_{j+1}^{-1}f^{(j)}).
\end{eqnarray*}
Hence
\begin{eqnarray*}
  \pi^{(N)}(e_0) (S_{j,j+1}-\hat{G}_{j,j+1})
  F_{\e_1\cdots\e_N}(z_1,\cdots,z_N) = \quad\qquad\qquad\qquad &&\\
  (S_{j,j+1}-\hat{G}_{j,j+1})
  (\sum_{k(\neq j,j+1)} e_0^{(k)}+ \hat{Y}_{j}^{-1}f^{(j+1)}+
  \hat{Y}_{j+1}^{-1}f^{(j)}) &&\\
  \times F_{\e_1\cdots\e_N}(z_1,\cdots,z_N) &\sim& 0,
\end{eqnarray*}
which completes the proof for $e_0$.

The case $x=f_0$ can be shown similarly by using
\[
S(t_1\otimes e_1)=(e_1\otimes 1)S.
\]

\subsection{Proof of \rf{rhof}}

The verification of \rf{rhof} is technically more complicated.
For this it is necessary to choose $p=q^4$.
In view of (\ref{eqn:Fcom}) it suffices to consider the case $j=N-1$.
The case $x=f_0$ being similar, we concentrate on the case $x=e_0$.

For convenience we use the decomposition of the $\ufop$-module
$V\otimes V=V^{(3)}\oplus V^{(1)}$,
the superscripts indicating the dimensions.
For $s=1$ or $3$ let
$v^{(s)}=\sum_{\e,\e'}v_\e\otimes v_{\e'}c_{\e\e'}^{(s)}$ be a vector in
$V^{(s)}$.
We have
\begin{equation}
Sv^{(3)}=-q^{-1}v^{(3)}, \qquad Sv^{(1)}=q v^{(1)}.
\label{eqn:Seig}
\end{equation}
Writing
$\sum_{\e,\e'}F_{\e_1\cdots\underbrace{\scriptstyle \e\e'}_{j,j+1}\cdots\e_N}
(z_1,\cdots,z_N)c_{\e\e'}^{(s)}$
as
$F_{\e_1\cdots c^{(s)}\cdots\e_N}(z_1,\cdots,z_N)$, we shall
verify \rf{rhof} separately for $s=1,3$.

\medskip

\noindent{\bf The case $s=3$:}\quad
{}From \rf{GENY} and \rf{SFREL}, we have
\begin{eqnarray*}
&&\rho^{(N)}(e_0)
F_{\e_1\cdots\e_N}(z_1,\cdots,z_N)\\
&&=\sum_{j=1}^Nf^{(N)}S_{N-1,N}\cdots S_{j,j+1}G^{-1}_{j-1,j}\cdots
G^{-1}_{1,2}
F_{\e_1\cdots\e_N}(z_2,\cdots,z_N,p^{-1}z_1).
\end{eqnarray*}
{}From (\ref{eqn:Ffus}) we have $F_{\e_1\cdots\e_{N-2}c^{(3)}}(z_1,\cdots,z_N)
\Big|_{z_N=q^{-2}z_{N-1}}\sim0$.
Therefore, we are to show
\begin{eqnarray*}
&&\sum_{j=1}^Nf^{(N)}S_{N-1,N}\cdots S_{j,j+1}G^{-1}_{j-1,j}\cdots G^{-1}_{1,2}
\\
&&\quad \times F_{\e_1\cdots \e_{N-2}c^{(3)}}(z_2,\cdots,z_N,p^{-1}z_1)
\res\sim0.
\end{eqnarray*}

{\it Step 1.}\quad
The sum of the terms for $j=N-1$ and $j=N$ reads
\bea
&&f^{(N)}(S_{N-1,N}+G^{-1}_{N-1,N})
G^{-1}_{N-2,N-1}\cdots G^{-1}_{1,2}\times\nonumber\\
&&\times F_{\e_1\cdots\e_{N-2}c^{(3)}}(z_2,\cdots,z_N,p^{-1}z_1)
\res\lb{M1}.
\ena
Because of (\ref{eqn:Seig}) and $G^{-1}_{N-1,N}\res=q^{-1}$,
the sum \rf{M1} vanishes.

{\it Step 2.}\quad
Note that
$S_{2,3}S_{1,2}(V\otimes V^{(3)})\subset V^{(3)}\otimes V$.
Therefore, for $j\le N-2$ the term
\be
f^{(N)}S_{N-1,N}\cdots S_{j,j+1}G^{-1}_{j-1,j}\cdots G^{-1}_{1,2}
F_{\e_1\cdots\e_{N-2}c^{(3)}}(z_2,\cdots,z_N,p^{-1}z_1)\lb{M2}
\en
is a linear combination of
\[
G^{-1}_{j-1,j}\cdots G^{-1}_{1,2}
F_{\e'_1\cdots\e'_{N-3}c'^{(3)}\e'_N}(z_2,\cdots,z_N,p^{-1}z_1)
\]
where $\sum_{\e,\e'}v_\e\otimes v_{\e'}c'^{(3)}_{\e\e'}\in V^{(3)}$.
Therefore, at $z_N=q^{-2}z_{N-1}$ the term \rf{M2} vanishes.

\medskip

\noindent{\bf The case $s=1$:}\quad
Next we consider the case of
$F_{\e_1\cdots\e_{N-2}c^{(1)}}$,
where
$v^{(1)}=\sum_{\e,\e'}v_\e\otimes v_{\e'}c^{(1)}_{\e\e'}$ is given by
(\ref{eqn:sing}).
Below, we give a complete proof for $N\le6$. We have verified
the general case by a similar method.

We set
\bea
B_{j,k}={q^{-1}z_j-qz_k\over z_j-z_k},\quad
C_{j,k}={(q-q^{-1})z_j\over z_j-z_k},\quad
{\overline C}_{j,k}={(q-q^{-1})z_k\over z_j-z_k},
\lb{BCC}
\ena
so that $G_{j,k}=B_{j,k}K_{j,k}+C_{j,k}$ and $G^{-1}_{j,k}=B_{j,k}K_{j,k}+
{\overline C}_{j,k}$.

\medskip
\begin{lem}\label{lem:1}
$S_{2,3}S_{1,2}v_\e\otimes v^{(1)}=q^{-1}v^{(1)}\otimes v_\e$.
\end{lem}

\begin{lem}\label{lem:2}
\begin{eqnarray*}
&&F_{\e_1\cdots\underbrace{\scriptstyle c^{(1)}}_{j,j+1}\cdots\e_N}
(z_1,\cdots,z_j,z_{j+1},\cdots,z_N)\Big|_{z_{j+1}=q^{-2}z_j}\\
&&=(-q)^{N-j}(1+q^{-2})
F_{\e_1\cdots,\e_{j-1},\e_{j+2}, \cdots\e_N}
(z_1,\cdots, z_{j-1},z_{j+2},\cdots,z_N)\times\\
&&\times\prod_{k=1}^{j-1}(z_k-q^2z_j)\prod_{k=j+2}^N(q^{-2}z_j-q^2z_k).
\end{eqnarray*}
\end{lem}

We are to show
\bea
&&\Bigl(q^2f^{(N)}\sum_{j=1}^N
S_{N-1,N}\cdots S_{j,j+1}G^{-1}_{j-1,j}\cdots G^{-1}_{1,2}
\nonumber\\
&&\times
F_{\e_1\cdots\e_{N-2}c^{(1)}}(z_2,\cdots,z_N,p^{-1}z_1)\Bigr)
\res\nonumber\\
&&\sim
-(q+q^{-1})f^{(N-2)}\sum_{j=1}^{N-2}
S_{N-3,N-2}\cdots S_{j,j+1}G^{-1}_{j-1,j}\cdots G^{-1}_{1,2}\nonumber\\
&&\times
\left(F_{\e_1\cdots\e_{N-2}}(z_2,\cdots,z_{N-2},p^{-1}z_1)
\prod_{k=1}^{N-2}(z_k-q^2z_{N-1})\right).
\lb{M4}
\ena

{\it Step 1.}\quad
By using Lemma \ref{lem:1} and Lemma \ref{lem:2},
the relation \rf{M4} reduces to
\bea
&&(1-q^2)z_{N-1}\sum_{j=1}^{N-2}
f^{(N-2)}S_{N-3,N-2}\cdots S_{j,j+1}G^{-1}_{j-1,j}\cdots G^{-1}_{1,2}
\times\nonumber\\
&&\times\Big\{F_{\e_1\cdots\e_{N-2}}(z_2,\cdots,z_{N-2},p^{-1}z_1)
\prod_{k=1}^{N-2}(z_k-q^2z_{N-1})\Big\}\nonumber\\
&&+q^2\Big\{
f^{(N)}G^{-1}_{N-2,N-1}\cdots G^{-1}_{1,2}\times\nonumber\\
&&F_{\e_1\cdots\e_{N-2}c^{(1)}}(z_2,\cdots,z_N,p^{-1}z_1)\Big\}\res
\sim0.\nonumber\\
\lb{M7}
\ena
If $N=2$, this is $0\sim0$. If $N=3$, the left hand side is
\be
(1-q^2)z_2f^{(1)}F_{\e_1}(p^{-1}z_1)
+q^2\Big\{f^{(3)}G^{-1}_{1,2}F_{\e_1c^{(1)}}(z_2,z_3,p^{-1}z_1)\Big\}
\Big|_{z_2=q^2z_3}.
\lb{M5}
\en
We have (see \rf{BCC})
\begin{eqnarray*}
&&f^{(3)}G^{-1}_{1,2}F_{\e_1c^{(1)}}(z_2,z_3,p^{-1}z_1)\\
&&=-q^{-1}(B_{1,2}K_{1,2}+\overline C_{1,2})
F_{\e_1-\,-}(z_2,z_3,p^{-1}z_1).
\end{eqnarray*}
We have $B_{1,2}K_{1,2}F_{\e_1-\,-}(z_2,z_3,p^{-1}z_1)
=B_{1,2}F_{\e_1-\,-}(z_1,z_3,p^{-1}z_2)$.
If $z_3=q^{-2}z_2$, then we have $p^{-1}z_2=q^{-2}z_3$. (Here we used $p=q^4$.)
Therefore, we have
\[
B_{1,2}K_{1,2}F_{\e_1-\,-}(z_2,z_3,p^{-1}z_1)\Big|_{z_3=q^{-2}z_2}\sim0.
\]

Noting that $F_{\e_1-\,-}(z_2,z_3,p^{-1}z_1)\Big|_{z_3=q^{-2}z_2}
=f^{(1)}F_{\e_1}(p^{-1}z_1)(z_2-z_1)$ and using \rf{BCC},
we obtain that \rf{M5} belongs to $\N$.

\medskip
\begin{lem}\label{lem:3}
${\overline C}_{2,3}G^{-1}_{1,2}=(G^{-1}_{1,2}+C_{2,3}){\overline C}_{1,3}$.
\end{lem}

\begin{lem}\label{lem:4}
\begin{eqnarray*}
&&B_{j,j+1}F_{\e_1\cdots\e_N}(z_1,\cdots,z_{j+1},z_j,\cdots,z_N)\\
&&\sim(S_{j,j+1}-C_{j,j+1})F_{\e_1\cdots\e_N}
(z_1,\cdots,z_j,z_{j+1},\cdots,z_N).
\end{eqnarray*}
\end{lem}

\noindent{\sl Remark.}\quad
In applying Lemma \ref{lem:4} more than once,
we must be careful about ordering the factors.
For example, for $N=5$ we have
\begin{eqnarray*}
&&B_{3,4}B_{2,4}F_{\e_1\e_2\e_3c^{(1)}}(z_4,z_2,z_3,z_5,p^{-1}z_1)\nonumber\\
&&=(S_{2,3}-C_{3,4})(S_{1,2}-C_{2,4})
F_{\e_1\e_2\e_3c^{(1)}}(z_2,z_3,z_4,z_5,p^{-1}z_1).
\end{eqnarray*}
\medskip

{\it Step 2.}\quad
By using Lemma \ref{lem:3} and Lemma \ref{lem:4}, we have
\bea
&&G^{-1}_{\scriptscriptstyle N-2,N-1}
\cdots G^{-1}_{\scriptscriptstyle 12}F_{\e_1\cdots\e_{N-2}c^{(1)}}
(z_2,\cdots,z_N,p^{-1}z_1)\res\nonumber\\
&&=\Big\{
(G^{-1}_{\scriptscriptstyle N-3,N-2}
+C_{\scriptscriptstyle N-2,N-1})(G^{-1}_{\scriptscriptstyle N-4,N-3}
+C_{\scriptscriptstyle N-3,N-1})\cdots(G^{-1}_{\scriptscriptstyle 12}
+C_{\scriptscriptstyle 2,N-1})\nonumber\\
&&+(G^{-1}_{\scriptscriptstyle N-4,N-3}+C_{\scriptscriptstyle N-3,N-1})
\cdots(G^{-1}_{\scriptscriptstyle 12}+C_{\scriptscriptstyle 2,N-1})
(S_{\scriptscriptstyle N-3,N-2}-C_{\scriptscriptstyle N-2,N-1})\nonumber\\
&&+(G^{-1}_{\scriptscriptstyle N-5,N-4}+C_{\scriptscriptstyle N-4,N-1})
\cdots(S_{\scriptscriptstyle N-3,N-2}-C_{\scriptscriptstyle N-2,N-1})
(S_{\scriptscriptstyle N-4,N-3}-C_{\scriptscriptstyle N-3,N-1})\nonumber\\
&&+\cdots\nonumber\\
&&+(S_{\scriptscriptstyle N-3,N-2}-C_{\scriptscriptstyle N-2,N-1})
(S_{\scriptscriptstyle N-4,N-3}-C_{\scriptscriptstyle N-3,N-1})
\cdots(S_{\scriptscriptstyle 12}-C_{\scriptscriptstyle 2,N-1})
\nonumber\Big\}\nonumber\\
&&\quad \times{\overline C}_{\scriptscriptstyle 1,N-1}
F_{\e_1\cdots\e_{N-2}c^{(1)}}(z_2,\cdots,z_N,p^{-1}z_1)\res
\nonumber\\
&&+
B_{N-2,N-1}\cdots B_{2,N-1}B_{1,N-1}
F_{\e_1,\cdots,\e_{N-2},c^{(1)}}
(z_1,\cdots,z_{N-2},z_N,p^{-1}z_{N-1})\Big|_{z_N=q^{-2}z_{N-1}}.
\nonumber\\
&&\lb{M6}
\ena
Note that the last term vanishes when $f^{(N)}$ is applied:
\[
f^{(N)}F_{\e_1,\cdots,\e_{N-2},c^{(1)}}
(z_1,\cdots,z_{N-2},z_N,p^{-1}z_{N-1})\Big|_{z_N=q^{-2}z_{N-1}}
\sim 0.
\]

For $N=5$, the derivation of \rf{M6} goes as follows.
By using Lemma \ref{lem:3} we have the equality
\begin{eqnarray*}&&B_{N-2,N-1}\cdots B_{2,N-1}B_{1,N-1}
G^{-1}_{3,4}G^{-1}_{2,3}G^{-1}_{1,2}\nonumber\\
&&={\overline C}_{3,4}G^{-1}_{2,3}G^{-1}_{1,2}
+B_{3,4}G^{-1}_{2,4}G^{-1}_{1,2}K_{3,4}
\nonumber\\
&&=(G^{-1}_{2,3}+C_{3,4})\oc_{2,4}G^{-1}_{1,2}
+B_{3,4}\oc_{2,4}G^{-1}_{1,2}K_{3,4}
+B_{3,4}B_{2,4}G^{-1}_{1,4}K_{2,4}K_{3,4}\nonumber\\
&&=(G^{-1}_{2,3}+C_{3,4})(G^{-1}_{1,2}+C_{2,4})\oc_{1,4}
+B_{3,4}(G^{-1}_{1,2}+C_{2,4})\oc_{1,4}K_{3,4}\nonumber\\
&&+B_{3,4}B_{2,4}\oc_{1,4}K_{2,4}K_{3,4}
+B_{3,4}B_{2,4}B_{1,4}K_{1,4}K_{2,4}K_{3,4}
\end{eqnarray*}
Using Lemma \ref{lem:4}, we have \rf{M6}.

{\it Step 3.}\quad
In the right hand side of \rf{M6},
if we pick up the terms
\begin{eqnarray*}
&&q^2f^{(N)}\Big\{\g32\g43\cdots G^{-1}_{2,3}G^{-1}_{1,2}\nonumber\\
&&+\g43\g54\cdots G^{-1}_{1,2}\es32\nonumber\\
&&+\g54\cdots G^{-1}_{1,2}\es32\es43\nonumber\\
&&+\cdots+\es32\es43\cdots S_{2,3}S_{1,2}\Big\}\times\nonumber\\
&&\times{\overline C}_{\scriptscriptstyle 1,N-1}
F_{\e_1\cdots\e_{N-2}c^{(1)}}(z_2,\cdots,z_N,p^{-1}z_1)\res,
\end{eqnarray*}
they cancel the first term in \rf{M7}.
For example, for $N=4$ we have
\begin{eqnarray*}
&&(1-q^2)z_3\Big\{f^{(2)}(S_{1,2}+G^{-1}_{1,2})F_{\e_1\e_2}(z_2,p^{-1}z_1)
(z_2-q^2z_3)\Big\}\nonumber\\
&&+q^2\Big\{f^{(4)}(G^{-1}_{1,2}+S_{1,2})\oc_{1,3}
F_{\e_1\e_2c^{(1)}}(z_2,z_3,z_4,p^{-1}z_1)
\Big\}\Big|_{z_4=q^{-2}z_3}\nonumber\\
&&\sim0.
\end{eqnarray*}

\begin{lem}\label{lem:5}
$(G^{-1}_{j,j+1}+C_{j+1,N-1})C_{j,N-1}=C_{j+1,N-1}G_{j,j+1}$.
\end{lem}

{\it Step 4.}\quad
Using Lemma \ref{lem:5},
we can show the cancellation of the terms obtained by setting
formally $G^{-1}_{1,2}=S_{N-3,N-2}=0$ in \rf{M6}.
For example, for $N=4$, we have $(C_{2,3}-C_{2,3})F_{\e_1\e_2c^{(1)}}
(z_2,z_3,z_4,p^{-1}z_1)\Big|_{z_4=q^{-2}z_3}=0$. For $N=5$, we have
\begin{eqnarray*}
&&\Big\{(G^{-1}_{2,3}+C_{3,4})C_{2,4}-C_{2,4}C_{3,4}-C_{3,4}
(S_{1,2}-C_{2,4})\Big\} \times\nonumber\\
&&\times F_{\e_1\e_2\e_3c^{(1)}}(z_2,z_3,z_4,z_5,p^{-1}z_1)
\Big|_{z_5=q^{-2}z_4}
\nonumber\\
&&\sim\Big\{C_{3,4}G_{2,3}-C_{2,4}C_{3,4}-C_{3,4}(G_{2,3}-C_{2,4})\Big\}\times
\nonumber\\
&&\times F_{\e_1\e_2\e_3c^{(1)}}(z_2,z_3,z_4,z_5,p^{-1}z_1)
\Big|_{z_5=q^{-2}z_4}\nonumber\\
&&=0.
\end{eqnarray*}
Here we used $S_{1,2}F_{\e_1\e_2\e_3c^{(1)}}(z_2,z_3,z_4,z_5,p^{-1}z_1)
\sim G_{2,3}F_{\e_1\e_2\e_3c^{(1)}}(z_2,z_3,z_4,z_5,p^{-1}z_1)$.

{\it Step 5.}\quad
The remaining terms on the right hand side of \rf{M6} are grouped into three
categories: the terms containing (i) $G^{-1}_{1,2}$ but not $S_{N-3,N-2}$,
(ii) $S_{N-3,N-2}$ but not $G^{-1}_{1,2}$,
(iii) both $G^{-1}_{1,2}$ and $S_{N-3,N-2}$.

The terms in (i), (ii), (iii) do not contain $K_{j,N-1}$, $K_{j,N}$
$(1\le j\le N-2)$. Therefore, if we apply these terms to
$F_{\e_1\cdots\e_{N-2}c^{(1)}}(z_2,\cdots,z_N,p^{-1}z_1)$ in \rf{M7},
the positions of $z_{N-1}$ and $z_N$ in
$F_{\e_1\cdots\e_{N-2}c^{(1)}}(z_2,\cdots,z_N,p^{-1}z_1)$ do not change.
Therefore, using
$\oc_{1,N-1}={\displaystyle{q^{-1}(1-q^2)z_{N-1}\over z_{N-1}-z_1}}$,
the relation \rf{M7} (to be proved) reduces to the form
\bea
&&\sum_{j=1}^{N-2}f^{(N-2)}\es32\cdots S_{j,j+1}G^{-1}_{j-1,j}\cdots
G^{-1}_{1,2}\Delta F\nonumber\\
&&-f^{(N-2)}\Big\{(\g32+C_{N-2,N-1})\cdots(G^{-1}_{1,2}+C_{2,N-1})\nonumber\\
&&+(\g43+C_{N-3,N-1})\cdots(\es32-C_{N-2,N-1})\nonumber\\
&&+\cdots+(S_{N-3,N-2}-C_{N-2,N-1})\cdots(S_{1,2}-C_{2,N-1})\Big\}
\Delta F\nonumber\\
&&\sim0,\lb{M8}
\ena
where we used the abbreviation
\begin{eqnarray*}
&&\Delta=\prod_{k=2}^{N-2}(z_k-q^2z_{N-1}),\nonumber\\
&&F=F_{\e_1\cdots\e_{N-2}}(z_2,\cdots,z_{N-2},p^{-1}z_1).
\end{eqnarray*}

Let us write $z_{1'}=p^{-1}z_1$ and
\be
G_{N-2,1'}={q^{-1}z_{N-2}-qz_{1'}\over z_{N-2}-z_{1'}}
K_{N-2,1}+{(q-q^{-1})z_{N-2}\over z_{N-2}-z_{1'}}.
\en
Then, the relation \rf{M8} reads
\begin{eqnarray*}
&&\sum_{j=1}^{N-2}f^{(N-2)}
\es32\cdots S_{j,j+1}G^{-1}_{j-1,j}\cdots G^{-1}_{1,2}\Delta F\nonumber\\
&&-f^{(N-2)}
\Big\{(\g32+C_{N-2,N-1})\cdots(G^{-1}_{1,2}+C_{2,N-1})\nonumber\\
&&+(\g43+C_{N-3,N-2})\cdots(G_{N-2,1'}-C_{N-2,N-1})\nonumber\\
&&+\cdots+(G_{N-2,1'}-C_{N-2,N-1})\cdots(G_{2,3}-C_{2,N-1})\Big\}\Delta F
\nonumber\\
&&\sim0.
\end{eqnarray*}

Let us call the terms in the first sum in the left hand side of this
relation ``the wanted terms'' or for short WT. The second sum contains WT.
The proof is over if we show the terms in the second sum other than WT
(we call them ``the unwanted terms''
or for short UWT), vanish. We will show this statement by induction on $N$.
For $N=2$ and $N=3$, the statement is trivial.
For $N=4$, the statement holds because $C_{2,3}-C_{2,3}=0$.
In general, we show the terms (i), (ii),(iii) other than WT separately cancel.
For example, for $N=5$ we have
\begin{displaymath}
(G^{-1}_{2,3}+C_{3,4})(G^{-1}_{1,2}+C_{2,4})
+(G^{-1}_{1,2}+C_{2,4})(G_{3,1'}-C_{3,4})
+(G_{3,1'}-C_{3,4})(G_{2,3}-C_{2,4}).
\end{displaymath}

UWT in (iii):0

UWT in (i)+(iii):
\begin{eqnarray*}
&&\Big\{(G^{-1}_{2,3}+C_{3,4})G^{-1}_{1,2}+G^{-1}_{1,2}(G_{3,1'}-C_{3,4})\Big\}
-(G^{-1}_{2,3}G^{-1}_{1,2}+G^{-1}_{1,2}G_{3,1'})\nonumber\\
&&=(C_{3,4}-C_{3,4})G^{-1}_{1,2}=0.
\end{eqnarray*}

UWT in (ii)+(iii):
\begin{eqnarray*}
&&\Big\{(G^{-1}_{1,2}+C_{2,4})G_{3,1'}+G_{3,1'}(G_{2,3}-C_{2,4})\Big\}
-(G^{-1}_{1,2}G_{3,1'}+G_{3,1'}G^{-1}_{2,3})\nonumber\\
&&=(C_{2,4}-C_{2,4})G_{3,1'}=0.
\end{eqnarray*}

For $N=6$ we have
\begin{eqnarray*}
&&(G^{-1}_{3,4}+C_{4,5})(G^{-1}_{2,3}+C_{3,5})(G^{-1}_{1,2}+C_{2,5})\nonumber\\
&&+(G^{-1}_{2,3}+C_{3,5})(G^{-1}_{1,2}+C_{2,5})(G_{4,1'}-C_{4,5})\nonumber\\
&&+(G^{-1}_{1,2}+C_{2,5})(G_{4,1'}-C_{4,5})(G_{3,4}-C_{3,5})\nonumber\\
&&+(G_{4,1'}-C_{4,5})(G_{3,4}-C_{3,5})(G_{2,3}-C_{2,5}).
\end{eqnarray*}

The cancellation of the terms in (iii) reduces to the case $N=4$.

The terms in (i)+(iii) are
\begin{eqnarray*}
&&(G^{-1}_{3,4}+C_{4,5})(G^{-1}_{2,3}+C_{3,5})G^{-1}_{1,2}\nonumber\\
&&+(G^{-1}_{2,3}+C_{3,5})G^{-1}_{1,2}(G_{4,1'}-C_{4,5})\nonumber\\
&&+G^{-1}_{1,2}(G_{4,1'}-C_{4,5})(G_{3,4}-C_{3,5}).\nonumber\\
\end{eqnarray*}
Except for the terms in (iii), which are already done,
we can commute $G^{-1}_{1,2}$ with the terms located on the right of
$G^{-1}_{1,2}$. By moving $G^{-1}_{1,2}$ to the right end,
the statement reduces to the cancellation of UWT in
\begin{displaymath}
  (G^{-1}_{3,4}+C_{4,5})(G^{-1}_{2,3}+C_{3,5})
+(G^{-1}_{2,3}+C_{3,5})(G_{4,1'}-C_{4,5})
+(G_{4,1'}-C_{4,5})(G_{3,4}-C_{3,5}).
\end{displaymath}
This is nothing but the case $N=5$. We omit the further details.

\section{Acknowledgements}
Two of us (R.~K. and T.~M.) thank the organizers of the meeting for
their hospitality. We are indebted to Akihiro Tsuchiya for his lecture
at Nagoya University on February 4, 1995. The work presented here is
an outcome of our efforts to understand his lecture. We also benefited
from the talk of Tomoki Nakanishi given at the International Institute
for Advanced Studies on November 21, 1994. We thank Katsuhisa Mimachi,
Atsushi Nakayashiki, Akihiro Tsuchiya and Yasuhiko Yamada for
discussions.  This work is partly supported by Grant-in-Aid for
Scientific Research on Priority Areas 231, the Ministry of Education,
Science and Culture.  H.~K. is supported by Soryushi Shyogakukai.
R.~K. and J.-U.~H.~P. are supported by the Japan Society for the
Promotion of Science.

\pagebreak[3]
\section{References}
\ifx\undefined\bysame
\newcommand{\bysame}{\leavevmode\hbox to3em{\hrulefill}\,}
\fi


\begin{thebibliography}{10}

\bibitem{Bkv}
A.~Berkovich, {\em Fermionic counting of {RSOS}-states and {Virasoro} character
  formulas for the unitary minimal series {M($\nu, \nu+1$)}. {Exact} results},
  Nucl. Phys. {\bf B431} (1994), 315.

\bibitem{BeMc}
A.~Berkovich and B.~M. McCoy, {\em Continued fractions and fermionic
  representations for characters of {$M(p,p')$} minimal models}, 1994, Preprint
  BONN-HE-94-28, ITPSB 94-060 hep-th/9412030.

\bibitem{BGHP}
D.~Bernard, M.~Gaudin, F.~D.~M. Haldane, and V.~Pasquier, {\em {Yang-Baxter}
  equation in long-range interacting systems}, J. Phys. A {\bf 26} (1993),
  5219--5236.

\bibitem{BPS}
D.~Bernard, V.~Pasquier, and D.~Serban, {\em Spinons in conformal field
  theory}, Nucl. Phys. {\bf B428} (1994), 612.

\bibitem{BLSb}
P.~Bouwknegt, A.~W.~W. Ludwig, and K.~Schoutens, {\em Spinon bases for higher
  level {SU(2)} {WZW} models}, 1994, Preprint USC-94/20, UCSB-TH-94, PUPT-1522,
  hep-th/9412108.

\bibitem{BLSa}
\bysame, {\em Spinon bases, {Yangian} symmetry and fermion representations of
  {Virasoro} characters in conformal field theory}, Phys. Lett. {\bf 338B}
  (1994), 448.

\bibitem{BLS94}
\bysame, {\em Spinon basis for $(\widehat{sl(2)})_k$ integrable highest weight
  modules}, preprint USC 94/22, PUPT-1537, 1994, to appear in the proceedings
  of `Statistical Mechanics and Quantum Field Theory,' USC, May 16-21, 1994.

\bibitem{Cal}
F.~Calogero, {\em Solution of the one-dimensional {$N$}-body problems
with quadratic
  and/or inversely quadratic pair potentials}, J. Math. Phys. {\bf 12} (1971),
  419--439.

\bibitem{ChPr90}
V.~Chari and A.~Pressley, {\em Yangians and {$R$}-matrices}, L'Enseigne\-ment
  Math{\'e}\-matique {\bf 36} (1990), 267--302.

\bibitem{CP}
\bysame, {\em Quantum affine algebras and affine {Hecke} algebras}, preprint,
  UCR and KCL, August 1993, q-alg/9501003.

\bibitem{Cher87}
I.~V. Cherednik, {\em A new interpretation of {Gelfand-Tzetlin} bases}, Duke.
  Math. J. {\bf 54} (1987), 563--577.

\bibitem{DFJMN}
B.~Davies, O.~Foda, M.~Jimbo, T.~Miwa, and A.~Nakayashiki, {\em Diagonalization
  of the {XXZ} {Hamiltonian} by vertex operators}, Comm. Math. Phys. {\bf 151}
  (1993), 89.

\bibitem{Dri85}
V.~G. Drinfeld, {\em Hopf algebras and the quantum {Yang-Baxter} equation},
  Soviet Math. Doklady {\bf 32} (1985), 254--258.

\bibitem{Dri86b}
\bysame, {\em Degenerate affine {Hecke} algebra and {Yangians}}, Funct. Anal.
  Appl. {\bf 20} (1986), 62--64.

\bibitem{Dri88}
\bysame, {\em A new realization of {Yangians} and quantized affine algebras},
  Soviet Math. Doklady {\bf 36} (1988), 212.

\bibitem{FaTa81}
L.~D. Faddeev and L.~A. Takhtadzhan, {\em What is the spin of a spin wave?},
  Phys. Lett. {\bf 85A} (1981), 375--377.

\bibitem{FaTa84}
\bysame, {\em Spectrum and scattering of excitations in the one-dimensional
  isotropic {Heisenberg} model}, J. Soviet Math. {\bf 24} (1984), 241--267.

\bibitem{FeON}
B.~L. Feigin, T.~Nakanishi, and H.~Ooguri, {\em The annihilating ideals of
  minimal models}, Proceedings of the RIMS Project 1991, Infinite Analysis
  (A.~Tsuchiya, T.~Eguchi, and M.~Jimbo, eds.), Advanced Series in Mathematical
  Physics, vol. 16A, World Scientific, Singapore, 1992, pp.~217--238.

\bibitem{hwm}
O.~Foda, K.~Iohara, M.~Jimbo, R.~Kedem, T.~Miwa, and H.~Yan, {\em Notes on
  highest weight modules of the elliptic algebra {${\cal
  A}_{q,p}(\widehat{sl}_2)$}}, 1994, RIMS preprint RIMS-979 (hep-th/9405058),
to appear in Proceedings 'Quantum Field Theory, Integrable Models
and Beyond', Eds. T. Inami and R. Sasaki, Suppl. Prog. Theor. Phys.

\bibitem{FoQu}
O.~Foda and Y-H. Quano, {\em Polynomial identities of the {Rogers} {Ramanujan}
  type}, 1994, Preprint, hep-th/9407191.

\bibitem{FR}
I.~B. Frenkel and N.~Yu Reshetikhin, {\em Quantum affine algebras and holonomic
  difference equations}, Comm. Math. Phys. {\bf 146} (1992), 1--60.

\bibitem{Hal}
F.~D.~M. Haldane, {\em Exact {Jastrow-Gutzwiller} resonating-valence-bond
  ground state of the spin-{$1/2$} antiferromagnetic {Heisenberg} chain with
  {$1/r^2$} exchange}, Phys. Rev. Lett. {\bf 60} (1988), 635--638.

\bibitem{HHTBP}
F.~D.~M. Haldane, Z.~N.~C. Ha, J.~C. Talstra, D.~Bernard, and V.~Pasquier, {\em
  Yangian symmetry of integrable quantum chains with long-range interactions
  and a new description of states in conformal field theory}, Phys. Rev. Lett
  {\bf 69} (1992), 2021--2025.

\bibitem{JM}
M.~Jimbo and T.~Miwa, {\em Algebraic analysis of solvable lattice models}, CBMS
  Regional Conference Series in Mathematics vol. {\bf 85}, AMS, 1994.

\bibitem{KKMM1}
R.~Kedem, T.~R. Klassen, B.~M. McCoy, and E.~Melzer, {\em Fermionic
  quasiparticle representations for characters of {$(G^{(1)})_1 \times
  (G^{(1)})_1 /(G^{(1)})_2$}}, Phys. Lett. {\bf 304} (1993), 263--270.

\bibitem{SUNY}
R.~Kedem and B.~M. McCoy, {\em Construction of modular branching functions from
  {Bethe's} equations in the 3-state {Potts} chain}, J. Stat. Phys. {\bf 71}
  (1993), 865.

\bibitem{KNS}
A.~Kuniba, T.~Nakanishi, and J.~Suzuki, {\em Characters in conformal field
  theories from thermodynamic {Bethe} ansatz}, Mod. Phys. Lett {\bf A8} (1993),
  1649.

\bibitem{LP}
J.~Lepowsky and M.~Primc, {\em Structure of the standard
modules for the affine Lie algebra $A^{(1)}_1$},
Contemporary Mathematics {\bf 46} AMS, Providence 1985.

\bibitem{Mel}
E.~Melzer, {\em The many faces of a character}, Lett. Math. Phys. {\bf 31}
  (1994), 233.

\bibitem{BONN}
W.~Nahm, A.~Recknagel, and M.~Terhoeven, {\em Dilogarithm identities in
  conformal field theory}, Mod. Phys. Lett. {\bf A8} (1993), 1835--1848.

\bibitem{NaYa}
A.~Nakayashiki and Y.~Yamada, {\em Crystalizing the spinon basis}, 1995,
  Preprint, hep-th/9504052.

\bibitem{Shas}
B.~S. Shastry, {\em Exact solution of an {$s=1/2$} {Heisenberg}
  antiferromagnetic chain with long-range interactions}, Phys. Rev. Lett. {\bf
  60} (1988), 639--642.

\bibitem{Suth}
B.~Sutherland, {\em Exact results for a quantum many--body problem
in one dimension. II}, Phys. Rev. A {\bf 5} (1972), 1372.

\bibitem{War}
S.~Ole Warnaar, {\em Fermionic solution of the {Andrews-Baxter-Forrester}
  model: unification of {TBA} and {CTM} methods}, 1995, Preprint,
  hep-th/9501134.

\end{thebibliography}
\end{document}